\documentclass[journal]{IAENGtran}

\usepackage{bm}
\usepackage{amsfonts, amssymb}

\usepackage{cite}

\ifCLASSINFOpdf
   \usepackage[pdftex]{graphicx}
   \DeclareGraphicsExtensions{.pdf,.jpeg,.png}
\else
   \usepackage[dvips]{graphicx}
   \DeclareGraphicsExtensions{.eps}
\fi

\usepackage[cmex10]{amsmath}
\usepackage{mathtools}

\begin{document}

\title{Systematic Improvement of Splitting Methods for the Hamilton Equations}

\author{Asif~Mushtaq, Anne Kv{\ae}rn{\o}, and K{\aa}re~Olaussen
 \thanks{Manuscript received April 16, 2012.}
 \thanks{A. Mushtaq is with the Department
 of Mathematical Sciences, NTNU, N-7048 Trondheim, Norway.
 e-mail: Asif.Mushtaq@math.ntnu.no.}
 \thanks{A. Kv{\ae}rn{\o} is with the Department
 of Mathematical Sciences, NTNU.
 e-mail: Anne@math.ntnu.no.}
 \thanks{K. Olaussen is with the Department
 of Physics, NTNU.
 e-mail: Kare Olaussen@ntnu.no.}}

\maketitle

\pagestyle{empty}
\thispagestyle{empty}

\begin{abstract}
We show how the standard (St{\"o}rmer-Verlet) splitting method for differential equations
of Hamiltonian mechanics (with accuracy of order $\tau^2$ 
for a timestep of length $\tau$) can be improved in a systematic
manner without using the composition method.  We give the explicit
expressions which increase the accuracy to order $\tau^8$, and
demonstrate that the method work on a simple anharmonic oscillator.

\end{abstract}

\begin{IAENGkeywords}
Splitting-method, Hamilton-equations, Higher-order-accuracy, Symplecticity
\end{IAENGkeywords}

\IAENGpeerreviewmaketitle

\section{Introduction}

\IAENGPARstart{T}{he} Hamilton equations of motion constitute a system of
ordinary first order differential equations,
\begin{align}
   \dot{q}^a = \frac{\partial H}{\partial p_a},\quad
   \dot{p}_a = -\frac{\partial H}{\partial q^a}, \quad a=1,\ldots,N,
   \label{HamiltonEquation}
\end{align}
where the $\dot{\ }$ denotes differentiation with respect to time $t$, and
$H = H(\bm{q}, \bm{p})$.  They can be viewed as the characteristic
equations of the partial differential equation
\begin{equation}
    \frac{\partial}{\partial t} \rho(\bm{q}, \bm{p}; t) =
    {\cal L}\,\rho(\bm{q}, \bm{p}; t),
    \label{HamiltonianFlow}
\end{equation}
with ${\cal L}$ the first order differential operator
\begin{equation}
    {\cal L} = \sum_{a=1}^N
    \frac{\partial H}{\partial p_a}\frac{\partial}{\partial q^a}
    -\frac{\partial H}{\partial q^a}\frac{\partial}{\partial p_a},
\end{equation}
generating a flow on phase space.
If $H$ does not depend ex\-pli\-citly on $t$, a formal solution
of (\ref{HamiltonianFlow}) is
\begin{equation}
    \rho(\bm{q}, \bm{p}; t) = 
    \text{e}^{t{\cal L}} \rho(\bm{q}, \bm{p}; 0).
\end{equation}
In most cases this expression remain just formal, but one may
often split the Hamiltonian into two parts,
$H = H_1 + H_2$, with a corresponding splitting
${\cal L} = {\cal L}_1 + {\cal L}_2$ such that the flows generated by
${\cal L}_1$ and ${\cal L}_2$ separately are integrable. One may
then use the Cambell-Baker-Hausdorff formula to approximate
the flow generated by ${\cal L}$. One obtains the Strang splitting formula ~\cite{HLWG,Strang}
\begin{equation}
     \text{e}^{\frac{1}{2}\tau{\cal L}_2}\,\text{e}^{\tau{\cal L}_1}\,
     \text{e}^{\frac{1}{2}\tau{\cal L}_2} = 
     \text{e}^{\tau{\cal L} + \frac{1}{24}\tau^3 \left[2{\cal L}_1+{\cal L}_2,\left[{\cal L}_1,{\cal L}_2\right]\right]+\cdots},
\end{equation}
which shows that time stepping this expression with a timestep $\tau$
provides an approximation with relative accuracy of order $\tau^2$,
exactly preserving the symplectic property of the flow.

This corresponds to the symplectic splitting scheme of iterating the process of solving
\begin{align}
       \dot{q}^a &= \frac{\partial H_2}{\partial p_a} ,\quad
       \dot{p}_a  = -\frac{\partial H_2}{\partial q^a},\quad\text{for
         step $\frac{1}{2}\tau$,}\nonumber\\
       \dot{q}^a &= \frac{\partial H_1}{\partial p_a} ,\quad
       \dot{p}_a  = -\frac{\partial H_1}{\partial q^a},\quad\text{for
         step $\tau$,}\label{SymplecticSplitting}\\
       \dot{q}^a &= \frac{\partial H_2}{\partial p_a} ,\quad
       \dot{p}_a  = -\frac{\partial H_2}{\partial q^a},\quad\text{for
         step $\frac{1}{2}\tau$.}\nonumber
\end{align}
Here the last part of one iteration may be combined with the
first part of the next, unless one deals with time
dependent systems or wants to register the state of
the system at the intermediate times.

From a practical point of view the most interesting property of
this formulation is that it can be interpreted directly
in terms of physical processes.
For instance, for Hamiltonians $H(\bm{q},\bm{p}) = T(\bm{p}) + V(\bm{q})$,
a standard splitting scheme is to choose $H_1 = T$ and $H_2 = V$. In that
case~(\ref{SymplecticSplitting}) corresponds to a collection of freely streaming
particles receiving kicks at regular time intervals $\tau$, these kicks being
dependent of the positions $\bm{q}$ of the particles. I.e, we may
think of the evolution as a collection of {\em kicks\/} and {\em moves\/} \cite{DTKC}.

It is not clear that this is the best way to approximate or model
the exact dynamics of the real system. 
For instance, why should the motion
between kicks be the free streaming generated by $T(\bm{p})$? 
There are more ways to split the Hamiltonian into two
integrable parts \cite{RIM}; the best splitting is most likely the one
which best mimics the physics of
equation~(\ref{HamiltonEquation}).
Further, since this equation is not solved exactly by (\ref{SymplecticSplitting})
for any finite value of $\tau$ we need not necessarily choose
$H_2$ to be {\em exactly\/} $H-H_1$ as long as it approaches this 
quantity sufficiently fast as $\tau\to 0$. We will
exploit this observation to improve the accuracy of the
splitting scheme (\ref{SymplecticSplitting}) in a
systematic manner. 

We are, of course, not the first trying to improve on the
St{\"o}rmer-Verlet splitting scheme. An accessible review of several
earlier approaches can be found in reference~\cite{BlanesCasasMurua}.
Neri~\cite{Neri}  has provided the general idea to construct
symplectic integrators for Hamiltonian systems.
Forest and Ruth ~\cite{FR} discussed the explicit fouth
order method for the integration of Hamiltonian equations for the
simplest non-trivial case. Yoshida~\cite{Yoshida} worked out a symplectic
integrator for any even order, and  Suzuki~\cite{Suzuki} presented the idea of
how recursive construction of successive approximants may be
extended to other methods.

\section{Harmonic Oscillators}

For a simple illustration of our idea consider the Hamiltonian
\begin{align}
   H({p},{q})=\frac{1}{2}\left( p^2 + q^2\right),
\end{align}
 whose exact evolution over a time interval $\tau$ is
\begin{equation}
     \begin{pmatrix} q_{\text{e}} \\[0.4ex]  p^{\text{e}} \end{pmatrix} =
    \begin{pmatrix*}[r] \cos \tau &  \sin \tau \\  - \sin \tau &
  \cos\tau \end{pmatrix*}   
   \begin{pmatrix} q\\p\end{pmatrix}.
\label{H_0a}
\end{equation}
Compare this with a {\em kick-move-kick\/} splitting 
scheme over the same time interval, with 
$H_{\text{kick}}= {\frac{1}{2} k q^2}$ and  $H_{\text{move}}=
{\frac{1}{2} m p^2}$, where $k$ and $m$ may depend on $\tau$.
One full iteration gives
\begin{equation}
    \begin{pmatrix}q_{\text{s}} \\ p^{\text{s}}\end{pmatrix} =
    \begin{pmatrix*}[c] 1-\frac{1}{2}m k\tau^2 & m \tau \\[0.4ex]
      -(1-\frac{1}{4} km \tau^ 2) k\tau &  1-\frac{1}{2} k m\tau^2 
    \end{pmatrix*}
    \begin{pmatrix}q \\ p \end{pmatrix}.
\end{equation}
We note that by choosing
\begin{align}
    m &= \frac{\sin\tau}{\tau} = 1 - \frac{1}{6}\tau^2 +
    \frac{1}{120}\tau^ 4 - \frac{1}{5040} \tau^6 + \cdots,\nonumber\\[-1ex]
    \label{HarmonicOscillatorCorrection}\\[-1ex]
    k &= \frac{2}{\tau} \tan \frac{\tau}{2} =
    1+\frac{1}{12}\tau^2+\frac{1}{120}\tau^4 + \frac{17}{20160}\tau^6
    + \cdots,\nonumber
\end{align}
the exact evolution is reproduced. If we instead choose a
{\em move-kick-move} splitting scheme, 
with $H_{\text{move}}=
{\frac{1}{2} \bar{m} p^2}$ and $H_{\text{kick}}= {\frac{1}{2} \bar{k}  q^2}$,
one iteration gives
\begin{equation}
    \begin{pmatrix}q_{\text{s}} \\[0.4ex] p^{\text{s}}\end{pmatrix} =
    \begin{pmatrix*}[c] 1-\frac{1}{2}\bar{m} \bar{k}\tau^2 & (1-\frac{1}{4}\bar{m}\bar{k})\bar{m} \tau \\[0.5ex]
      -\bar{k}\tau &  1-\frac{1}{2} \bar{k}\bar{m}\tau^2 
    \end{pmatrix*}
    \begin{pmatrix}q \\ p \end{pmatrix},
\end{equation}
which becomes exact if we choose
\begin{equation}
  \bar{m} = \frac{2}{\tau}\tan\frac{\tau}{2},\quad
  \bar{k} = \frac{\sin\tau}{\tau}.
\end{equation}
It should be clear that this idea works for systems of harmonic
oscillators in general, i.e. for quadratic Hamiltonians of the form
\begin{equation}
   H(\bm{q},\bm{p}) = \frac{1}{2} \left( \bm{p}^T M \bm{p} + \bm{q}^T
     K \bm{q} \right),
\end{equation}
where $M$ and $K$ are symmetric matrices.
For a choosen splitting scheme and step 
interval $\tau$ there are always modified matrices
$M_\tau = M + {\cal O}(\tau^2)$ and $K_\tau = K + {\cal O}(\tau^2)$
which reproduces the exact time evolution. For systems where
$M$ and $K$ are too large for exact diagonalization,
but sparse,  a systematic expansion of $M_\tau$ and $K_\tau$
in powers of $\tau^2$ could be an efficient way to improve the
standard splitting schemes.

\section{Nonlinear systems}

For a more general treatment we consider Hamiltonians of the form
\begin{equation}
    H(\bm{q}, \bm{p}) = \frac{1}{2} \bm{p}^T M \bm{p}  + V(\bm{q}).
\end{equation}
A series solution of the Hamilton equations in powers of $\tau$ is
\begin{align}
     q^a_{\text{e}} &= q^a + p^a\tau - \frac{1}{2}\partial^a V \tau^2 
     -\frac{1}{6}\partial^a (DV) \tau^3+{\cal O}(\tau^4),\nonumber\\
     p^{\text{e}}_a &= p_a - \partial_a V \tau  -
     \frac{1}{2}\partial_a (DV) \tau^2 \\
     &+\partial_a \left( \frac{1}{12} \bar{D}V -\frac{1}{6} D^2V\right) \tau^3 + {\cal O}(\tau^4)\nonumber
\end{align}
Here we have introduced notation to shorten expressions,
\begin{align}
     &\partial_a \equiv \frac{\partial}{\partial q^a},\quad 
     \partial^a \equiv M^{ab}\partial_b,\quad
     p^a \equiv M^{ab} p_b,
      \nonumber\\[-1.5ex]
\\[-1.5ex]      
     &D \equiv p_a \partial^a,\quad
     \bar{D} \equiv (\partial_a V)\partial^a,\nonumber
\end{align}
where we employ the {\em Einstein summation convention\/}: An
index which occur twice, once in lower position and once in upper
position, are implicitly summed over all available values. I.e, 
$M^{ab} \partial_b \equiv \sum_b M^{ab} \partial_b$ (we will generally
use the matrix $M$ to rise an index from lower to upper position).
The corresponding result for the {\em kick-move-kick\/}
splitting scheme is
\begin{align}
    q_{\text{s}}^a &= q^a + p^a\tau - \frac{1}{2}\partial^a V \tau^2  +{\cal O}(\tau^4),\nonumber\\
    p^{\text{s}}_a &=p_a - \partial_a V\tau -
     \frac{1}{2}\partial_a (DV) \tau^2 \\
     &+\partial_a \left( \frac{1}{8}\bar{D}V  - \frac{1}{4}D^2V\right) \tau^3 + {\cal O}(\tau^4).\nonumber
\end{align}
As expected it differs from the exact result in the third order, but
the difference can be corrected by introducing second order
generators
\begin{equation}
    T_2 = -\frac{1}{12} D^2V \tau^2,\quad
    V_2 = \frac{1}{24} \bar{D}V\tau^2,
\end{equation}
to be used in respectively the {\em move\/} and {\em kick} steps.
Specialized to a one-dimensional system with potential $V = \frac{1}{2}q^2$
this agrees with equation~(\ref{HarmonicOscillatorCorrection}). With
this correction the {\em kick-move-kick\/} splitting scheme agrees
with the exact solution to $4^{\text{th}}$ order in $\tau$, but differ
in the $\tau^5$-terms. We may correct the difference by introducing
fourth order generators,
\begin{align}
    T_4 &=\frac{1}{720}\left( D^4  - 9 \bar{D} D^2 + 3 D\bar{D} D\right) V\tau^4,\nonumber\\[-1ex]
  \\[-1ex]
    V_4 &=\frac{1}{480} \bar{D}^2 V\tau^4.\nonumber
\end{align}
Specialized to a one-dimensional system with potential $V = \frac{1}{2}q^2$
this agrees with equation~(\ref{HarmonicOscillatorCorrection}). 
With this correction the {\em kick-move-kick\/} splitting scheme agrees
with the exact solution to $6^{\text{th}}$ order in $\tau$, but differ
in the $\tau^7$-terms. We may correct the difference by introducing
sixth order generators,
\begin{align}
    T_6 &= -\frac{1}{60480}\left(2\, D^6  \right.
         -40\, \bar{D}D^4
         +46\, D\bar{D} D^3\nonumber\\ &
        -15\,D^2 \bar{D}D^2 +54\,\bar{D}^2 D^2 - 9\,\bar{D}D\bar{D}D\nonumber\\ 
       &- 42\,D\bar{D}^2 D\left. +12\, D^2\bar{D}^2\right)V\tau^6\nonumber\\[-1ex]
  \\[-1ex]
    V_6 &= \frac{1}{161280}\left( 17\, \bar{D}^3 - 10\,\bar{D}_3\right)V\tau^6,\nonumber
\end{align}
where we have introduced 
\begin{equation}
 \bar{D}_3 \equiv (\partial_aV)(\partial_b V)(\partial_c
 V) \partial^a \partial^b \partial^c.
\end{equation}
Specialized to a one-dimensional system with potential $V = \frac{1}{2}q^2$
this agrees with equation~(\ref{HarmonicOscillatorCorrection}). 
With this correction the {\em kick-move-kick\/} splitting scheme agrees
with the exact solution to $8^{\text{th}}$ order in $\tau$, but differ
in the $\tau^9$-terms. One may continue the correction process, but
this is probably well beyond the limit of practical use already.

\section{Solving the {\em move\/} steps}

Addition of extra potential terms $V \rightarrow V_{\text{eff}} \equiv V + V_2 + V_4 + \dots$ 
is in principle unproblematic for solution of the {\em kick\/} steps. The
equations,
\begin{equation}
  \dot{q}^a = 0,\quad
  \dot{p}_a = -\partial_a V_{\text{eff}}(\bm{q}),
\end{equation}
can still be integrated exactly, preserving the symplectic structure.
The situation is different for the kinectic term
$T \rightarrow T_{\text{eff}} \equiv T + T_2 + T_4 + \cdots$,
since it now leads to equations
\begin{equation}
  \dot{q}^a = \frac{\partial }{\partial p_a}T_{\text{eff}}(\bm{q},\bm{p}),\quad
  \dot{p}_a = -\partial_a T_{\text{eff}}(\bm{q},\bm{p}),
  \label{MoveSteps}
\end{equation}
which is no longer straightforward to integrate exactly.
Although the problematic terms are small
one should make sure that the {\em move\/} steps preserve
the symplectic structure exactly. Let $\bm{q}, \bm{p}$
denote the positions and momenta just before the {\em move\/} step, and $\bm{Q}, \bm{P}$
the positions and momenta just after. We construct a generating 
function~\cite{Goldstein, Arnold, LubichEtAlVI51}
$G(\bm{q}, \bm{P};\tau)$, with
\begin{equation}
   Q^a = \frac{\partial G}{\partial P_a},\quad
   p_a = \frac{\partial G}{\partial q^a}.
   \label{CanonicalTransformation}
\end{equation}
This preserves the symplectic structure; we just have to construct
$G$ to represent the {\em move\/} step sufficiently accurately.
Consider first the case without the correction terms. The choice
$G = q^a P_a + \frac{1}{2} P^a P_a \tau$ gives
\begin{equation}
     Q^a = q^a + P_a\,\tau,\quad p_a = P_a,
     \label{SimpleKick}
\end{equation}
which is the correct relation. Now add the $T_2$-term to 
the {\em move\/} step. To order $\tau^4$ the exact solution
of equation (\ref{MoveSteps}) becomes
\begin{align}
    Q^a &= q^a + p_a\,\tau - \frac{1}{6} \partial^a DV\, \tau^3 -\frac{1}{24}\partial^a D^2V\, \tau^4,\nonumber\\[-1.5ex]
    \label{ExactMove} \\[-1.5ex]
    P_a &= p_a + \frac{1}{12} \partial_a \, D^2 V\,\tau^3 + \frac{1}{24}\partial_a D^3 V\,\tau^4.\nonumber
\end{align}
Compare this with the result of changing
\begin{equation}
  G \rightarrow G -\frac{1}{12} {\cal D}^2\, V \tau^3 - \frac{1}{24} {\cal D}^3V\, \tau^4,
\end{equation}
where ${\cal D} \equiv P_a \partial^a$. The solution of equation~(\ref{CanonicalTransformation})
change from the relations~(\ref{SimpleKick}) to
\begin{align}
  Q^a &= q^a + P^a\tau - \frac{1}{6}\partial^a {\cal D}V \tau^3 
  -\frac{1}{8}\partial^a {\cal D}^2V\,\tau^4,\label{Qequation}\\
  p_a &= P_a - \frac{1}{12} \partial_a {\cal D}^2 V \tau^3 - \frac{1}{24}\partial_a{\cal D}^3 V \tau^4. \label{Pequation}
\end{align}
Since ${\cal D}$ is  linear in $\bm{P}$,
equation~(\ref{Pequation}) constitute a system of third order
algebraic equation which in general must be solved numerically.
This should usually be a fast process for small $\tau$.
An exact solution of this equation is required to preserve
the symplectic structure, but  this solution should also
agree with the exact solution of (\ref{MoveSteps}) to
order $\tau^4$. This may be verified by perturbation
expansion in $\tau$. A perturbative solution of equation~(\ref{Pequation}) is
\begin{equation*}
  P_a = p_a + \frac{1}{12} \partial_a D^2 V \tau^3 + \frac{1}{24} \partial_a D^3 V \tau^4 + \ldots,
\end{equation*}
which inserted into (\ref{Qequation}) reproduces the full solution~(\ref{ExactMove}) to order $\tau^4$.

This process can be systematically continued to higher orders.
We write the transformation function as
\begin{equation}
    G(\tau) = \sum_{n=0}^\infty G_n\,\tau^n,
\end{equation}
and find the first terms in the expansion to be
{\footnotesize
\begin{align}
    G_0 &= q^a P_a,\nonumber\\
    G_1 &= \frac{1}{2} P^a P_a,\nonumber\\
    G_2 &= 0,\nonumber\\
    G_3 &= -\frac{1}{12}{\cal D}^2 V,\nonumber\\ 
    G_4 &= -\frac{1}{24}{\cal D}^3 V,\nonumber\\
    G_5 &= -\frac{1}{240}\left( 3\,{\cal D}^4+ 3\, \bar{D}{\cal D}^2
      -{\cal D} \bar{D} {\cal D} \right) V,\\
    G_6 &= -\frac{1}{720} \left({2\,\cal D}^5  + 8\,\bar{D}{\cal D}^3
      - 5\, {\cal D}\bar{D}{\cal D}^2 \right) V,\nonumber\\
    G_7 &=-\frac{1}{20160}\left({10\,\cal D}^6 + 10\,\bar{D}{\cal D}^4 + 90\,{\cal
        D}\bar{D} {\cal D}^3 - 75\,{\cal D}^2 \bar{D} {\cal D}^2\right.
    \nonumber
    \\ &\phantom{=-\frac{1}{20160}}\left.\, +18\,\bar{D}^2{\cal D}^2 -3\, \bar{D}{\cal
         D}\bar{D}{\cal D} -14\,{\cal D}\bar{D}^2{\cal D} + 4\,
       {\cal D}^2\bar{D}^2 \right) V,\nonumber\\
    G_8 &= - \frac{1}{40320} \left(3\,{\cal D}^7
      -87\,\bar{D}{\cal D}^5
    +231\,{\cal D}\bar{D}{\cal D}^4
    -133\,{\cal D}^2\bar {D}{\cal D}^3 
    \right.\nonumber\\
  &\phantom{= - \frac{1}{40320}} \left.\, +63\, \bar{D}^2 {\cal D}^3
    -3\,{\cal D}\bar{D}^2{\cal D}^2
     -21\,{\cal D}^2 \bar{D}^2 {\cal D}
      +4\,{\cal D}^3 \bar{D}^2
      \right.\nonumber\\ &\phantom{= - \frac{1}{40320}}\left.\,-63
    \,\bar{D}{\cal D}\bar{D}{\cal D}^2
       +25\,{\cal D}\bar{D}{\cal D}\bar{D}{\cal D}\right) V.\nonumber
\end{align}
}

\section{Explict computations}

\begin{figure}[!h]
\begin{center}
\includegraphics[clip, trim = 8ex 6ex 9ex 5ex, width=0.483\textwidth]{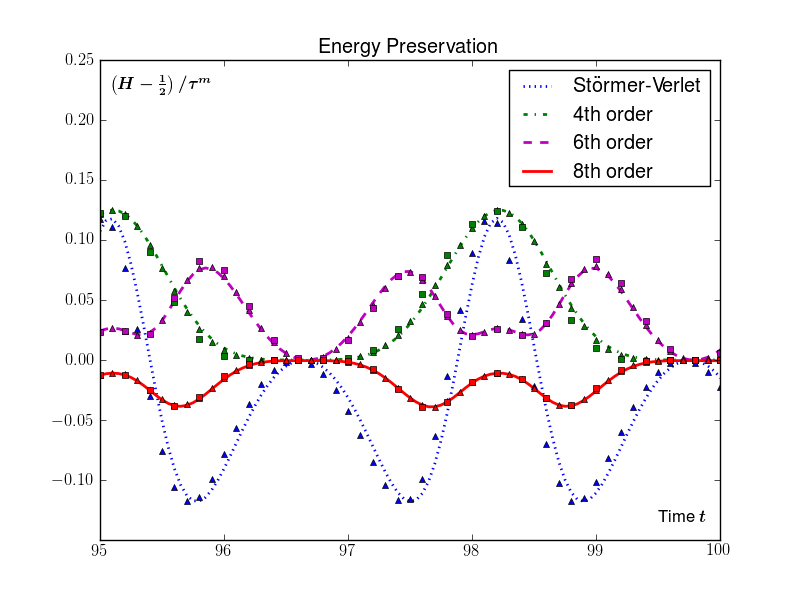}
\end{center}
\caption{This figure illustrate how well energy is conserved with the
various splitting schemes. The quanties plotted is $(H-\frac{1}{2})/\tau^{m}$
for $\tau = 0.2$ (squares), $\tau=0.1$ (triangles) and $\tau=0.05$ (lines).
Here $m=2$ for the St{\"o}rmer-Verlet scheme (dotted line), $m=4$ for the
$\tau^2$-corrected generators (dash-dotted line), $m=6$ for the
$\tau^4$-corrected generators (dashed line), and $m=8$  for the
$\tau^6$-corrected generators (fulldrawn line). Each plotted quantity is essentially the
value of the next correction at the visited point in phase space. Since the plot is taken
over the last half of the $16^{\text{th}}$ period the figure also give some indication of how
well the exact oscillation period is reproduced by the scheme. The deviation is quite
large for the St{\"o}rmer-Verlet scheme when $\tau=0.2$; to avoid cluttering the figure
we have not included these points.
}
\label{energyPreservation}
\end{figure}

It remains to demonstrate
that our algorithms can be applied to real examples. We have
considered the Hamiltonian
\begin{equation}
   H = \frac{1}{2} p^2 + \frac{1}{4}q^4,
\end{equation}
with initial condition $q(0)=0$, $p(0)=1$. The exact motion is a
nonlinear oscillation with $H$ constant equal to $\frac{1}{2}$,
and period  
\begin{equation}
    T = 4\,\int_0^{2^{1/4}} \frac{\sqrt{2}\, \text{d}q}{\sqrt{2-q^4}} = 2^{1/4}\,
\text{B}(\frac{1}{4},\frac{1}{2}) \approx 6.236\,339\ldots.
\end{equation}
Here $\text{B}(x,y) = \Gamma(x)\Gamma(y)/\Gamma(x+y)$ is the beta function.
In figure~\ref{energyPreservation} we plot the behaviour
of $\left(H-\frac{1}{2}\right)/\tau^{2+n}$ during the last half of the $16^{\text{th}}$ oscillation,
for various values of $\tau$ and corrected generators up to order $\tau^6$ (corresponding to $n=6$).

\begin{figure}[!h]
\begin{center}
\includegraphics[clip, trim = 8ex 6ex 9ex 5ex, width=0.483\textwidth]{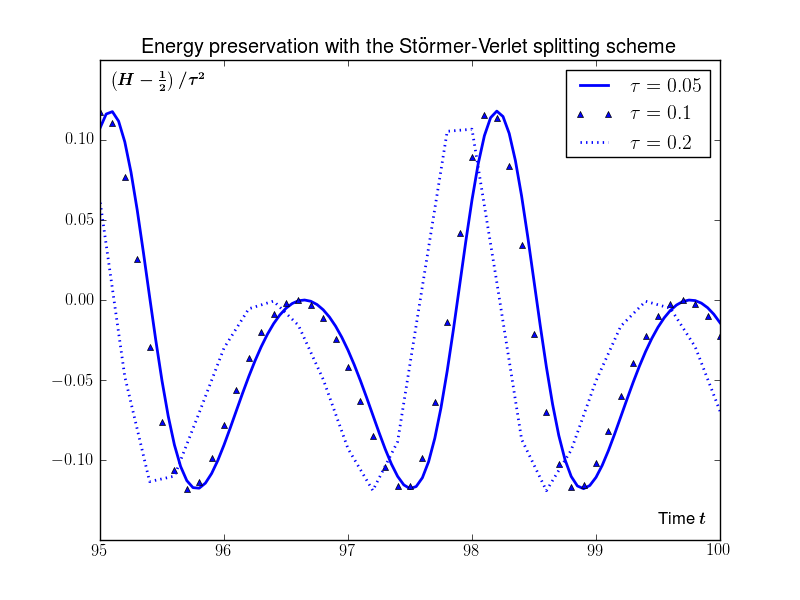}
\end{center}
\caption{This figure illustrate the long time behaviour (through the
  last half of the $16^{\text{th}}$ period) for
the St{\"o}rmer-Verlet scheme. Different timesteps $\tau$
have an effect on the period of oscillation, but the preservation of energy
remains stable for a very long time.
}
\label{stormerVerlet}
\end{figure}

\begin{figure}[!h]
\begin{center}
\includegraphics[clip, trim = 8ex 6ex 9ex 5ex, width=0.483\textwidth]{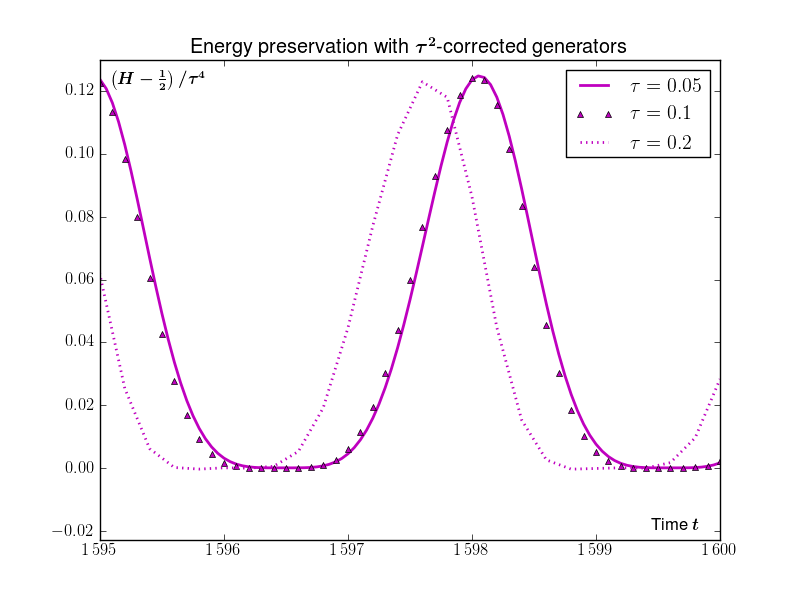}
\end{center}
\caption{This figure illustrate the long time behaviour (through the first
  half of the $257^{\text{th}}$ period) for
the $\tau^2$-corrected scheme. Different timesteps $\tau$
have an effect on the period of oscillation, but the preservation of energy
remains stable for a very long time.
}
\label{fourthOrder}
\end{figure}

\begin{figure}[!h]
\begin{center}
\includegraphics[clip, trim = 8ex 6ex 9ex 5ex, width=0.483\textwidth]{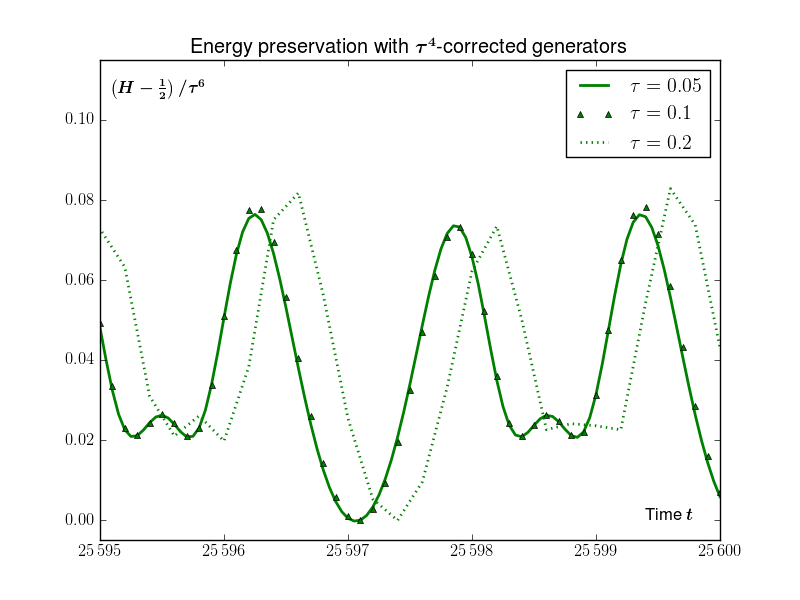}
\end{center}
\caption{
This figure illustrate the long time behaviour for
the $\tau^4$-corrected scheme (through the last half of the $4\,104^{\text{th}}$ period). Different timesteps $\tau$
have an effect on the period of oscillation, but the preservation of energy
remains stable for a very long time. 
}
\label{sixthOrder}
\end{figure}

\begin{figure}[!h]
\begin{center}
\includegraphics[clip, trim = 8ex 6ex 9ex 5ex, width=0.483\textwidth]{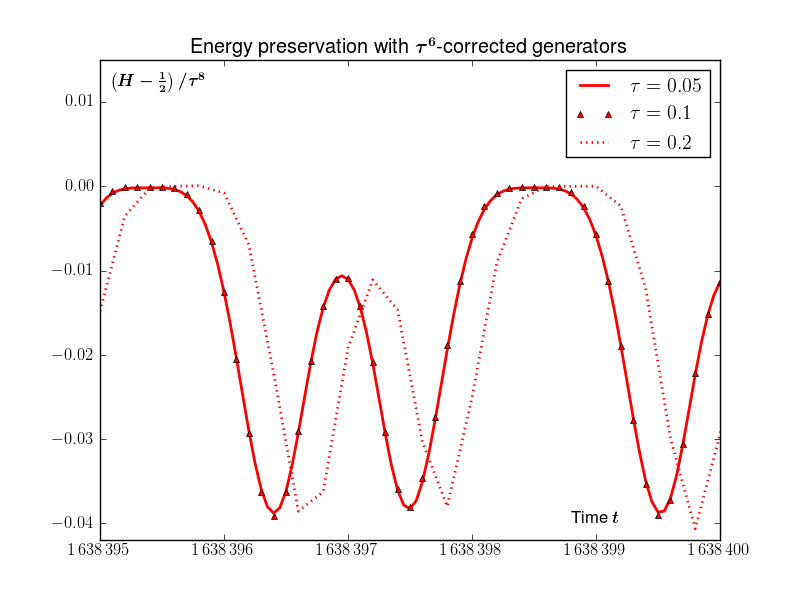}
\end{center}
\caption{
This figure illustrate the long time behaviour for
the $\tau^6$-corrected scheme (through the last half of period $262\,718$). Different timesteps $\tau$
have an effect on the period of oscillation, but the preservation of energy
remains stable for a very long time (for high accuracy and very long runs
the effect of numerical roundoff errors eventually becomes visible). 
}
\label{eightOrder}
\end{figure}


\section{Conclusion}

We have shown that it is possible to systematically improve
the accuracy of the usual symplectic integration schemes for
a rather general class of Hamilton equations. The process
is quite simple for linear equations, where it may be useful
for sparse systems. For general systems the method requires
the solution of a set of nonlinear algebraic equations
at each {\em move\/} step. To which extent an higher-order
method is advantageous or not will depend on the system
under analysis, and the wanted accuracy. As always
with higher order methods the increased accuracy per
step may be countered by the higher computational cost
per step \cite{AAK}.

\ifCLASSOPTIONcaptionsoff
  \newpage
\fi

\newpage


%

\end{document}